\documentclass{article} 

\usepackage{authblk} \usepackage{amsmath, amsfonts,
	mathpazo} \usepackage{geometry}

\usepackage{cite} \usepackage{graphicx} \usepackage{amsmath}
\usepackage{algorithm} 
\usepackage{algpseudocode} \usepackage{graphics} \usepackage{epsfig}
\usepackage{comment} \usepackage{subfigure}
\usepackage{url}
\usepackage{tikz} \usetikzlibrary{positioning} 
\usepackage{xcolor} 
\usepackage{pstricks}

\definecolor{mycolor}{rgb}{0.1,0.5,1}

\newtheorem{remark}{Remark} 

\usepackage{booktabs}
\newcommand\keywords[1]{\textbf{Keywords}: #1}

\usepackage{multicol}
\numberwithin{figure}{section}
\numberwithin{table}{section}

\usepackage{listings}
\usepackage{authblk}
\title{An MP-DWR method for $h$-adaptive finite element methods}

\author[1]{Chengyu Liu} \author[1,2,3]{Guanghui Hu}

\affil[1]{Department of Mathematics, Faculty of Science and
	Technology, University of Macau, Macao SAR, China.}
\affil[2]{Zhuhai UM Science \& Technology Research Institute, Zhuhai,
	Guangdong, China.}
\affil[3]{Guangdong-Hong Kong-Macao Joint Laboratory for Data-Driven
	Fluid Mechanics and Engineering Applications, University of Macau,
	Macao SAR, China.}

\begin{document}

\maketitle
\abstract{In a dual weighted residual method based on the finite
	element framework, the Galerkin orthogonality is an issue
	that prevents solving the dual equation in the same space as
	the one for the primal equation. In the literature, there have
	been two popular approaches to constructing a new space for the
	dual problem, i.e., refining mesh grids ($h$-approach) and
	raising the order of approximate polynomials
	($p$-approach). In this paper, a novel approach is proposed
	for the purpose based on the multiple-precision technique,
	i.e., the construction of the new finite element space is
	based on the same configuration as the one for the primal
	equation, except for the precision in calculations. The feasibility of
	such a new approach is discussed in detail in the paper. In numerical experiments, the proposed approach can be realized conveniently with C++ \textit{template}. Moreover, the new approach shows remarkable improvements in both efficiency and storage compared with the $h$-approach and the
	$p$-approach. It is worth mentioning that the performance of our approach is comparable
	with the one through a higher order interpolation
	($i$-approach) in the literature. The combination of these two
	approaches is believed to further enhance the
	efficiency of the dual weighted residual method.}

\section{Introduction}\label{sec1}
The use of adjoint equations and duality arguments in residual-based \textit{a posteriori} error estimation has become one of the most popular topics in the numerical analysis of partial differential equations. The idea can be traced back to the work of Babu\v{s}ka and Miller \cite{babuvska1984post,babuvska1984post2,babuvska1987feedback}. Then the techniques were systematically developed for more general situations \cite{eriksson1988adaptive,eriksson1991adaptive,eriksson1995introduction}. Furthermore, Becker and Rannacher improved this approach into a computation-based feedback method, i.e. the dual weighted residual (DWR) method, which has the goal-orientation property \cite{DWRfirst,becker1996weighted,becker2001optimal,rannacher2001adaptive}. Related techniques using duality arguments in post-processing and design have been studied in \cite{peraire1998bounds,prudhomme1999goal,giles1997adaptive,giles2002adjoint}.
Due to the capability of dealing with the quantities of interest, the DWR method has been applied in many areas, such as mechanics, physics, and chemistry \cite{DWRBook,Rannacher1999, bangerth2010adaptive,kormann_2016,bruchhauser2017numerical,rabizadeh2020pointwise,avijit2022efficient}.

Based on a variational formulation of the discrete problem, the DWR method uses duality techniques to generate \textit{a posteriori} error estimation, which depends on the unknown exact solution of a so-called dual problem. It is noted that the dual solution can not be calculated in the same finite element space as the primal problem. Otherwise, the error estimation would be vanished due to the Galerkin orthogonality. To deal with this problem, there are some classic approaches to obtain an approximation of the dual solution, such as finer mesh approximation ($h$-approach), higher-order finite element method approximation ($p$-approach), or patch-wise higher-order interpolation approximation ($i$-approach)\cite{DWRBook,giles2002adjoint}. In these classic approaches, several operations are necessary to construct the dual solution, which would cost unignorable computational resources. To overcome this, in this work, we explore the possibility of building a different finite element space for solving the dual problem through one feature of the computer, i.e., finite precision arithmetic. 

\begin{figure}[H]
  \begin{center}
    \fbox{
      \includegraphics[width=0.5\textwidth]{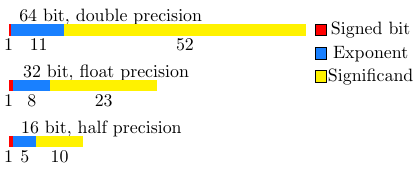}	
    }	
    \caption{IEEE754 data format.}
      \label{IEEE754}
  \end{center}
\end{figure}

A key observation is that with the same configuration (same mesh, degrees of freedom, basis functions, etc.), if we use two different floating point precision to build finite element spaces, for example, float-precision and double-precision, the two finite element spaces should be different from the computer's point of view. To understand it, we need to recall the way a real number is expressed in the computer. From IEEE standard shown in Fig. \ref{IEEE754}, each real number takes up 64 bits in double-precision formats, 32 bits in float-precision formats, and 16 bits in half-precision, respectively. Consequently, one number stored with different precision formats will be treated as different numbers in calculations.
Then we take two vectors $\boldsymbol{v}_1=(cos\frac{\pi}{3},sin\frac{\pi}{3})$
and $\boldsymbol{v}_2=(-sin\frac{\pi}{3},cos\frac{\pi}{3})$ as an example, it can be seen that two vectors are orthogonal with each other
mathematically. Numerically, if two vectors are expressed in single-precision,
the inner product $(\boldsymbol{v}_1,\boldsymbol{v}_2)$ is computed as 0.00000000. 
However, if $\boldsymbol{v}_1$
is in single, and $\boldsymbol{v}_2$ is in double, then the inner product becomes
-0.000000007771812, which is obviously not zero in double-precision. Such an
observation motivates an idea for the new implementation of DWR, i.e., solving
primal and dual problems in two different precision formats to avoid the Galerkin
orthogonality.

In this paper, a multiple-precision DWR (MP-DWR) method is proposed for the
implementation of $h$-adaptive finite element methods, towards improving the
simulation efficiency. In the framework of the MP-DWR method, the primal and dual problem are solved with different precision, respectively. Based on these, a DWR \textit{a posteriori}
error estimation can be generated to serve the following local mesh refinement. The feasibility of above algorithm is discussed in detail in the paper, from which the break of the Galerkin orthogonality can be seen clearly in examples.

The improvement for the efficiency can be expected from our method by observations that i). both the re-partitioning of mesh grids and the rearrangement of degree of freedoms as well as associate basis functions are avoided in our method, 
and ii). the involvement of single-precision calculation would bring the reduction of both CPU time and storage. 
In numerical experiments, the significant acceleration by our method can be seen clearly, e.g., in solving 3D Poisson equation, around an order of magnitude difference of CPU time between the $h$-approach strategy and our method can be observed in the generation of the error indicator.

Another feature of our MP-DWR method is its concise coding work with the aid of \textit{template} property in C++ language. In this paper, all simulations are
realized by using a C++ library \textit{AFEPack}\cite{li2005multi}, in which the code is organized following a standard procedure on solving partial differential equations as modules. By using these modules and setting the type of the floating-point number as a parameter of the template class, the code can be concise and easily managed. 

The outline of this paper is as follows. In Section 2, based on the Poisson
equation, the idea of goal-orientation and the dual weighted residual method in
$h$-adaptive finite element methods are briefly reviewed. In Section 3, we
describe a novel DWR implementation based on multiple-precision in detail. In
Section 4, several numerical experiments are implemented, in which the validity and acceleration of the proposed method can be demonstrated. 
Besides that, the capability of the proposed method for solving the eigenvalue problem is also shown. In addition, several remarks are delivered about the properties and
limitations of the proposed method, based on our numerical experience. Finally, we give the conclusion of this paper and future works in Section 5.

\section{Poisson equation and DWR based $h$-adaptive mesh methods}\label{sec2}

In this section, based on the Poisson equation, a DWR based $h$-adaptive finite
element method is reviewed briefly, as well as the associated classic theory.

\subsection{Poisson equation and finite element discretization}
Consider the Poisson equation 

	\begin{equation}
		\label{Poisson}
		\left\{
		\begin{array}{rl}
			-\Delta u = f,&\text{ in } \Omega,\\
			u = 0,&\text{ on } \partial\Omega,
		\end{array}
		\right.
	\end{equation}
 
with a homogeneous Dirichlet boundary condition on the boundary
$\partial \Omega$ for the well posedness. $\Omega\subset \mathbb{R}^d$
is a bounded domain in two ($d = 2$) or three ($d = 3$) dimensional space. We choose
the variational space $V= H_{0}^{1}(\Omega)$ to incorporate the homogeneous
Dirichlet boundary condition. This implies $f\in L^{2}(\Omega)$.

Then the weak solution $u$ to the Poisson equation can be characterized via the variational problem:
 
\begin{equation}
	Find\ u\in V\ ,\ such\ that\ a(u,v):=(\nabla u,\nabla v) = (f,v),\ \forall v\in V. 
	\label{weak form}
\end{equation}
  
On a triangulation $\mathcal{T}$, we can introduce the finite element space
$V_h\subset H^{1}_{0}(\Omega)$. 
In $V_h$, the Ritz-Galerkin approximation to (\ref{weak form}) can be expressed as:
 
\begin{equation}
	\label{GalerkinDiscretization}
	Find\ u_h\in V_h,\ such\ that\ a(u_h,v):=(\nabla u_h,\nabla v)=(f,v),\ \forall v\in V_h\subset V,
\end{equation}
 
where $u_h$ is the Galerkin approximation to $u$.

\subsection{Goal-orientation and dual weighted residual}

We now introduce the classic residual type \textit{a posteriori} error
estimation. Based on that, an introduction of the framework of a
goal-oriented \textit{a posteriori} error estimation based on DWR method
follows. Denote the exact solution of (\ref{Poisson}) and (\ref{GalerkinDiscretization})
by $u$ and $u_h$, respectively. Then the error $e_h$ is
related to the residual as:
 
	\begin{equation}
		\begin{array}{l}
			\displaystyle a(e_h,v)=\displaystyle a(u-u_h,v):=\mathcal{R}(v),\quad\forall v\in V.
		\end{array}
		\label{a(eh,v)}
	\end{equation}
 
Then integration by parts element-wise to transform $\mathcal{R}$ as:
 
\begin{equation}
	\label{R}
	\mathcal{R}(v)=\sum_{K}\int_K(f+\Delta u_h)vd\boldsymbol{x}+\sum_{E}\int_E r[\boldsymbol{n}_E\cdot\nabla u_h]vds,
\end{equation}
 
where $K$ and $E$ represent elements and the edges of elements,
respectively, $\boldsymbol{n}_E$ represents a unit vector on $E$, and $r$ is the ``jump" term as defined in \cite{verfurth2013posteriori}.
Based on (\ref{R}), a classic residual-based error indicator $\eta_{K}^{\textit{resi}}$ for every element
$K\in \ \mathcal{T}$ can be obtained as given in \cite{verfurth2013posteriori}.

Generally, the quantities of interest can be some functionals of the exact solution $u$ in simulations. The dual weighted residual method is a better choice to deal with these situations, in which the residual terms are multiplied by weights obtained by solving a dual problem. Following by \cite{rannacher2001adaptive,DWRBook}, the framework of the DWR method can be briefly summarized as follows. Based on the above finite element discretization, to evaluate the target functional $J(\cdot)$ (assumed to be linear for simplicity), the corresponding dual problem is introduced as
 
\begin{equation}
	\label{dual}
	a(\varphi,z)=J(\varphi),\ \ \forall \varphi\in V,
\end{equation}
 
where $z\in V$ is the corresponding dual solution. The Galerkin approximation of the dual problem in $V_h$ can be given by 
 
\begin{equation}
	\label{Galerkin}
	a(\varphi_h,z_h)=J(\varphi_h),\ \ \forall \varphi_h\in V_h,
\end{equation}
 
where $z_h\in V_h$. By combining (\ref{dual}) and (\ref{Galerkin}), we can get the dual weighted residual error indicator of error in functional $J(\cdot)$:
 
\begin{equation}
		\label{dwrRepre}
	\begin{array}{rl}
		J(e):&=J(u) - J(u_h) = a(u-u_h,z) \\
		&= \sum\limits_{K\in \mathcal{T}_h}\{(R_h,z)_{K}+\frac{1}{2}(r(\boldsymbol{n}\cdot\nabla u_h), z)_{\partial K}\} \vspace{1ex} =:\eta, 
	\end{array}	
\end{equation}
 
where the cell residual $R_h$ is defined by
 
\begin{equation}
	R_{h\vert K} := f+\Delta u_h,
\end{equation}
 
and $r(\boldsymbol{n}\cdot\nabla u_h)$ represents the jump term of neighboring elements $K$ and $K'$ on the common edge $\Gamma$ with normal unit vector $\boldsymbol{n}$ pointing from $K$ to $K'$:
 
	\begin{equation}
		r(\boldsymbol{n}\cdot\nabla u_h) = (\nabla u_{h\vert K'\cap \Gamma}-\nabla u_{h\vert K\cap\Gamma})\cdot \boldsymbol{n}.
	\end{equation}
 
In simulations, with the approximation $\tilde{\eta}$ to $\eta$, the effectivity of the error indicator can be evaluated through the index given by
 
\begin{equation}
	\label{EffInd}
	I_{\textit{eff}} = \left\vert \frac{\tilde{\eta}}{J(e)}\right\vert,
\end{equation}
which represents the degree of overestimation and should desirably be close to one \cite{DWRBook}.
 
Moreover, \textit{a posterori} error estimate can be obtained as
 
\begin{equation}
	\label{dwrError}
	\vert J(e) \vert \leq \eta_{\omega} := \sum_{K\in \mathcal{T}_h} \rho_K \omega_K,
\end{equation}
 
where the cell residuals $\rho_K$ and weights $\omega_K$ are given by:
 
\begin{equation}
	\begin{array}{rl}
		&\rho_{K} := (\Vert R_h\Vert_{K}^{2}+h_{K}^{-1}\Vert r(\boldsymbol{n}\cdot\nabla u_h) \Vert^{2}_{\partial K})^{1/2}, \vspace{1ex} \\
		&\omega_{K} := (\Vert z \Vert_{K}^{2}+h_K \Vert z\Vert^{2}_{\partial K})^{1/2}.
	\end{array}	
\end{equation}

Since the dual solution $z$ is generally unknown, the only left issue is to obtain the approximation $\tilde{z}$. It should be noticed that $\tilde{z}$ can not be computed in the same space $V_h$ as the one for solving the primal problem otherwise the error indicator would be a trivial value due to Galerkin orthogonality:
 
\begin{equation}
	\label{GO}
	a(u-u_h,v_h) = 0,\ \forall v_h \in V_h.
\end{equation}

To deal with this issue, there are three main classic approaches for computing $\tilde{z}$ as follows.

\begin{itemize}
	\item By $h$-approach: A denser mesh $\tilde{\mathcal{T}_h}$ can be constructed by
	globally refining the mesh $\mathcal{T}_h$. Based on $\tilde{\mathcal{T}_h}$, a
	different space $\tilde{V_h}$ can be built. Then the dual solution $\tilde{z}$ can be
	obtained from $\tilde{V_h}$.
	
	\item By $p$-approach: A finite element space $\tilde{V}_{h}^{(2)}$ with quadratic (or higher-order) finite elements has to be constructed. From $\tilde{V}_{h}^{(2)}$, the dual
	solution $\tilde{z}$ can be obtained.
	
	\item By $i$-approach: The bilinear Ritz projection $z_h \in V_h$ of $z$ is computed first. Then the approximation $\tilde{z}$ is obtained by patch-wise higher-order interpolation of $z_h$. In the following numerical experiments, we only consider \textit{quadratic} interpolation.
\end{itemize}

In practical simulations, the first two approaches all build a different
finite element space for solving the dual problem. Some operations would be necessary in such a process, i.e., re-meshing the grid, re-organizing the degree of freedoms, regeneration of the information of numerical quadrature,the reformation of a larger scale stiff matrix, and solving an even larger linear system of linear equations, 
which would cost nonignorable CPU time and memory storage in
calculations. To avoid these operations, the $i$-approach obtains the approximation $\tilde{z}$ by using higher-order interpolation, which makes it cheaper than the other two classic approaches \cite{rabizadeh2020pointwise}. Different from three classic approaches above, we expect to propose a novel approach to avoid the Galerkin orthogonality. Such a new approach can also avoid these operations so the acceleration in simulations can be expected.

\section{A novel DWR implementation based on multiple-precision}\label{sec3}
Recall the example in the introduction, which shows that the orthogonality of two vectors can not be preserved numerically with multiple-precision calculations. Based on the same idea, we introduce a new implementation of DWR, in which the multiple-precision will be used to effectively break the Galerkin orthogonality.

First of all, we will show how the multiple-precision breaks the Galerkin orthogonality to make sure this idea works.
The Galerkin orthogonality in (\ref{GO}) indicates that: for $v_h=u_h$ from $V_h$, we have 
 
\begin{equation}
	\label{GOnew}
	a(u-u_h, u_h) = 0.
\end{equation}
 
By using a different precision format, a different space $V_h'$ can be built without changing the mesh grid. From $V_h'$, $u_h'$ can be obtained. Based on the above idea, $u_h$ and $u_h'$ would be different since they are computed with different precision formats. Furthermore, they are obtained from two different finite element spaces $V_h$ and $V_h'$, respectively. Thus, it can be expected that the Galerkin orthogonality between them would not be preserved:
 
\begin{equation}
	\label{notGO}
	a(u-u_{h},u_{h}')\neq 0.
\end{equation}

To confirm that, we implement the numerical experiment by using a combination of double- and float-precision. Here, we only focus on the Poisson equation in the 2-dimensional domain $\Omega=(-1,1)^2$ with the exact function as $u = sin(\pi x)sin(2\pi y)$. With double-precision, the finite element space $V_h^d$ is built firstly, from which $u_h^d$ can be obtained. Then the similar operations are implemented in float precision, we can obtain the numerical solution $u_h^f$ from $V_h^f$. Based on these, we can compute the value of $a(u-u_h^d,u_h^d)$ and $a(u-u_h^d,u_h^f)$. The numerical results are shown in Table \ref{GO_check1}. 
\begin{table}[h]
  \begin{center}
    \begin{tabular}{@{}llll@{}}
      \toprule
      Dofs  	 & L2 Error  & $a(u-u_h^d,u_h^d)$ & $a(u-u_h^d,u_h^f)$\\ 
      \midrule
      385   & 4.470e-1   &  -5.035e-5      & -8.117e-5  \\
      1473  & 1.136e-1   & -1.579e-7       & -5.103e-4  \\
      5761  & 2.854e-2   & -6.596e-10      & -2.202e-3  \\
      22785 & 7.143e-3   & 8.254e-11       & -9.337e-3  \\
      \bottomrule
    \end{tabular}
  \end{center}
    \caption{Galerkin orthogonality of $u_h^d$ and $u_h^f$ for the case $u = sin(\pi x)sin(2\pi y)$.}
   	\label{GO_check1}%
\end{table}

Tables \ref{GO_check1} shows that the error in L2-norm (L2 Error) tends to zero with more Dofs (degrees of freedom), such that $u_h^d$ and $u_h^f$ get more accurate. Meanwhile, the value of $a(u-u_h^d,u_h^d)$ tends to zero, which means the Galerkin orthogonality of $u_h^d$ is numerically preserved better. That is the reason that we can not directly compute the dual solution in the same space with the same precision. On the contrary, the value of $a(u-u_h^d,u_h^f)$ gets further away from zero as the number of Dofs increases, which means the Galerkin orthogonality between $u_h^d$ and $u_h^f$ is numerically broken more obviously. Based on this phenomenon, we can confirm that the Galerkin orthogonality between $u_h^d$ and $u_h^f$ can be numerically broken by multiple-precision.


Thus, it can be expected that dual solutions from the finite element space constructed with different precision will work in the DWR method. By combining such primal and dual solutions, a reliable DWR error indicator can be obtained for mesh adaptation. Based on this, we propose Algorithm \ref{alg:Framwork} for the multiple-precision implementation of the DWR method.
\begin{algorithm}[htb]
  \caption{Framework of Adaptive finite element method based on MP-DWR(DF) method.} 
  \label{alg:Framwork} 
  \begin{algorithmic}[1] 
    \Require 
    The initial mesh, $\mathcal{T}_0$; 
    The tolerance, $tol$; 
    The maximum iteration times, $Max_{\text{iter}}$; the iteration times $k=0$.
    
    \State Build a finite element space $V_h$ on the mesh $\mathcal{T}_{k}$; 
    \State Solve the primal problem in double-precision to get the primal solution $u_{h}^d$; 
    \State Solve the dual problem in float-precision to get the dual solution $z_{h}^f$;	
    \State Calculate the cell residual $\tilde{\rho}_K$ and weights $\tilde{\omega}_K$ for each element to construct the error indicator $\tilde{\eta}_{\omega}$ as (\ref{dwrError});
    \State Check the stopping condition:
    If $\tilde{\eta}<tol$ or $k>Max_{\text{iter}}$, stop the adaptation and return the numerical solutions.
    \State Else, adapt the current mesh with $\tilde{\eta}_{\omega}$ to get the new mesh $\mathcal{T}_{k+1}$, update the iteration times $k=k+1$ and return to step 1.
  \end{algorithmic} 
  \label{Algorithm}
\end{algorithm}
\begin{remark}
	In the above framework of MP-DWR(\textbf{DF}) method, the primal problem is solved in \textbf{D}ouble-precision. As for dual problem, it is solved with \textbf{F}loat-precision. In fact, the framework of MP-DWR(FD) is also useful, i.e. solving the primal and dual problem in Float and Double precision, respectively. But there are some limits of MP-DWR(FD), which will be shown through a numerical example (\ref{limitofFD}) in the following part. In this paper, MP-DWR(DF) framework is used by default in numerical experiments if not specified.
\end{remark}

Such a new implementation of DWR method has a series of advantages, which can be summarized as follows.
First of all, several repeated operations of building finite element space can be avoided. Compared to the two classic methods ($h$-approach and $p$-approach), the MP-DWR method only computes the modules of building the finite element space once in one iteration. With the help of
multiple-precision, the re-partitioning of the mesh grid, the rearrangement of
degrees of freedom, the recalculation of the quadrature information, the
reformation of stiff matrix, and solving an even larger linear system of linear equations can be avoided in the process of obtaining the dual solution. In this way, considerable improvements in efficiency can be expected. Compared with the $i$-approach, the dual solution in MP-DWR is generally not as accurate as the one from higher-order interpolation, making the error indicator potentially less accurate. But the MP-DWR method can avoid high-order interpolation, thus saving some computation time.

The second advantage of MP-DWR comes from the involvement of computations with
low precision. On computers, if one number is stored in double-precision, it
will cost 8 bytes. While it costs 4 bytes in single-precision. Thus, if data are stored in single-precision instead of
double-precision, the storage memory would be saved by fifty percent. Besides that, data with lower precision are more conducive to accelerating the calculation. Consider following C++ codes:
\lstset{basicstyle=\small\ttfamily,frame=single,framexleftmargin=-1pt,framexrightmargin=0pt,framesep=1pt,linewidth=0.98\textwidth,language=c++}
\begin{lstlisting}
double double_data = 0.;
for(int i=0;i<1.0e7;i++){double_data+=4.0*atan(1.0);}
float float_data = 0.;
for(int i=0;i<1.0e7;i++){float_data+=float(4.0*atan(1.0));}
\end{lstlisting}
The same operation of the data is implemented 1.0e7 times with
different precision formats. For double-precision, the cost of CPU time is 4.010e-2. As a comparison, CPU time is 2.448e-2 for float-precision, which is much less than double-precision. Therefore, it can be expected that the single-precision calculation part in the MP-DWR method can bring further reduction in both storage and CPU time.  

In addition, the realization of the MP-DWR method can be very concise by using
\emph{template}. In this paper, we implement all simulations based on the
\textit{AFEPack} library \cite{li2005multi}. The code in the \textit{AFEPack} library follows a
standard procedure for solving partial differential equations. Furthermore, each
step is realized
as a module based on \textit{template} in \textit{AFEPack}. \textit{Template} is a feature of the C++ programming language that allows functions and classes to work on many different data types without rewriting the code for each one. By setting the floating type as the parameter of \textit{template}, it is quite
trivial to realize Algorithm \ref{Algorithm}. 
For example, consider the bilinear operator module in \textit{AFEPack}: $\text{BilinearOperator}<\dots,\text{Number}>$.
We can merely change the \emph{template} parameter `Number', then the precision of the numbers in \textit{BilinearOperator} will be changed. Thus, it is convenient to implement multiple-precision calculations without repeating the entire code for each precision format. Moreover, the use of \textit{template} also facilitates the modification of our codes based on some function libraries. Codes of one numerical example of the MP-DWR method can be found in Github \cite{GithubExample}.

\section{Numerical experiments}\label{sec4}
In this section, a series of examples are shown to demonstrate the validity and acceleration of the MP-DWR method. In the following numerical experiments, we use the linear triangle template
element for 2D cases and the tetrahedron template element for
3D cases, respectively. All numerical
examples presented in this paper are developed by using the
\textit{AFEPack}\cite{li2005multi} in C++, and the hardware is a DELL Precision T5610
workstation with 12 CPU Cores, 3.6GHz and 64 Gigabytes memory.

Moreover, a high order Gaussian quadrature is recommended in the calculations since it can produce more accurate results, which can break the Galerkin orthogonality more obviously. Therefore, the MP-DWR error indicators can be more effective. Here, we use the six-order Gaussian quadrature in all computations.

\subsection{Effectiveness of MP-DWR method}
\counterwithin{figure}{subsection}
\counterwithin{table}{subsection}
In this part, we conduct a series of numerical experiments based on MP-DWR method to verify its effectiveness. In the beginning, the Poisson equation is considered for testing the goal-orientation property and effectivity of the MP-DWR method. 
Then experiments of a problem with diffusion coefficient and the convection-diffusion equation are also conducted to study whether the MP-DWR method can be effective in some more complex situations.
\subsubsection{Poisson equation}\label{PoissonCase}
In this example, we consider the Poisson problem as (\ref{Poisson}) with $\Omega = (-1,1)\times(-1,1)$. First of all, we study the effectivity indices of the MP-DWR method for different target functionals in the Poisson case. As a benchmark example in \cite{DWRBook}, we consider an exact solution $u_{1}=(1-x_1^2)(1-x_2^2)sin(4x_1)sin(4x_2)$ with the target functionals:
 
\begin{equation}
	\label{Target4Poisson}
	\begin{array}{rl}
			J_{\text{area}}(u)&:=\vert S\vert^{-1}\int_S u d\boldsymbol{x},\ S:=[-\frac{1}{2},0]\times[0,\frac{1}{2}],\vspace{1ex} \ \boldsymbol{x}=(x_1,x_2),\\
			J_p(u)&:=u(\boldsymbol{x}_e),\ \boldsymbol{x}_e = (\frac{1}{2}, \frac{1}{2}).
	\end{array}
\end{equation}
 
\begin{figure}[h]
	\centering
	\includegraphics[width=0.8\textwidth]{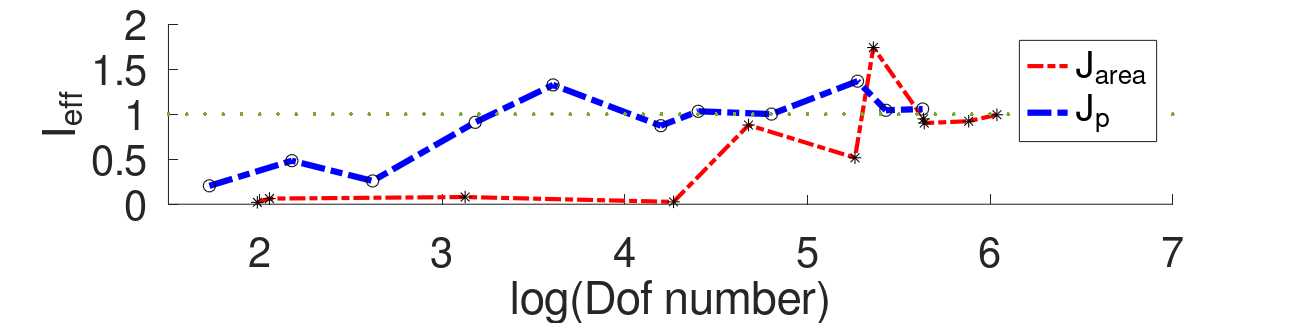}
	\caption{Effectivity indices of MP-DWR error indicators for the target functionals $J_{\text{area}}$ and $J_p$ over the number of Dofs.}
	\label{EffectIndices}
\end{figure}

To evaluate the error indicator generated by MP-DWR, we compute the effectivity index $I_{\text{eff}}$ given in (\ref{EffInd}). In Figure \ref{EffectIndices}, we present effect indices of the MP-DWR approach applied to this Poisson case. It can be seen that effectivity indices of target functionals $J_1$ and $J_p$ get closer to one as the number of Dofs increases. Although they are not strictly equal to one, they tend to around one. Therefore, we can confirm the effectivity of the MP-DWR method.

Then, a sharper function is considered, of which the exact form is given as $u_2=5.0\times(1-x_1^2)(1-x_2^2)exp(1-x_2^{-4})$. Besides that, we consider a point-wise output functional, which is given by:
 
\begin{equation}
	J_p(u):=u(\boldsymbol{x}_e),\ \boldsymbol{x}_e=(0.2,0.8).
\end{equation}
 
\begin{figure}
\begin{minipage}[c]{0.3\linewidth}
		\includegraphics[width=0.6\textwidth]{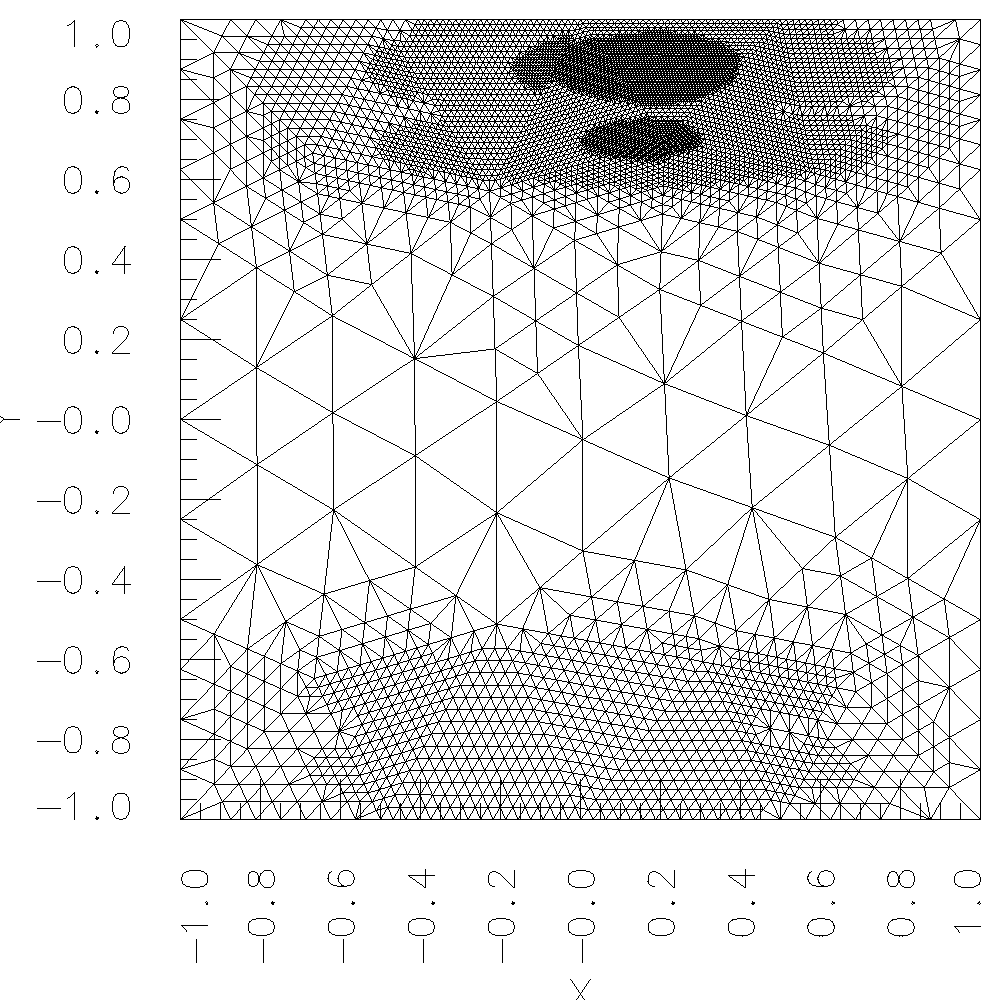}\\
		\includegraphics[width=0.6\textwidth]{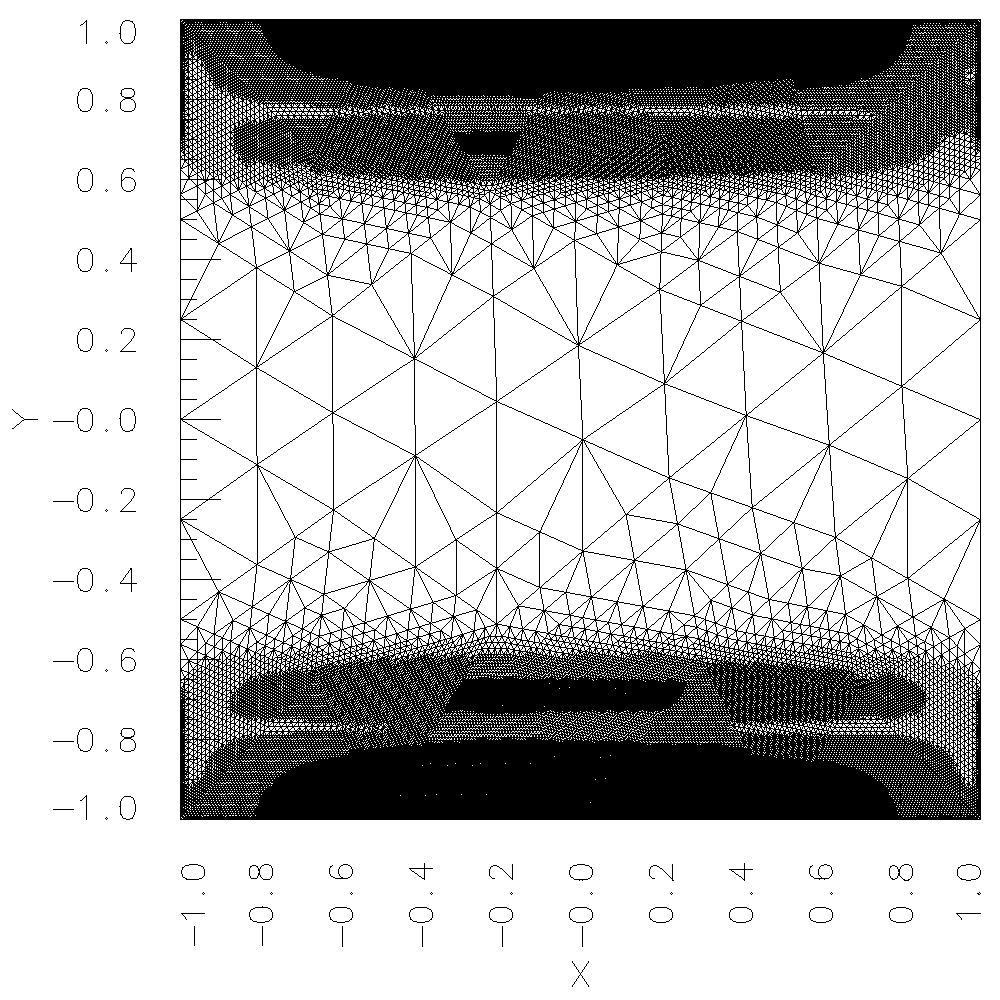}
\end{minipage}
\begin{minipage}[c]{0.7\linewidth}
		\includegraphics[width=0.5\textwidth]{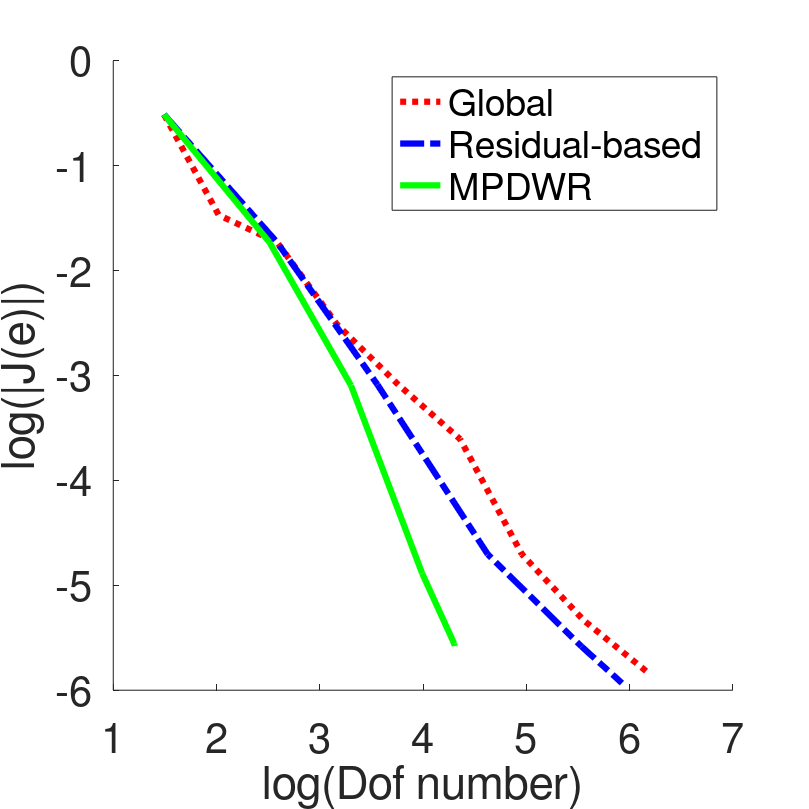}\includegraphics[width=0.5\textwidth]{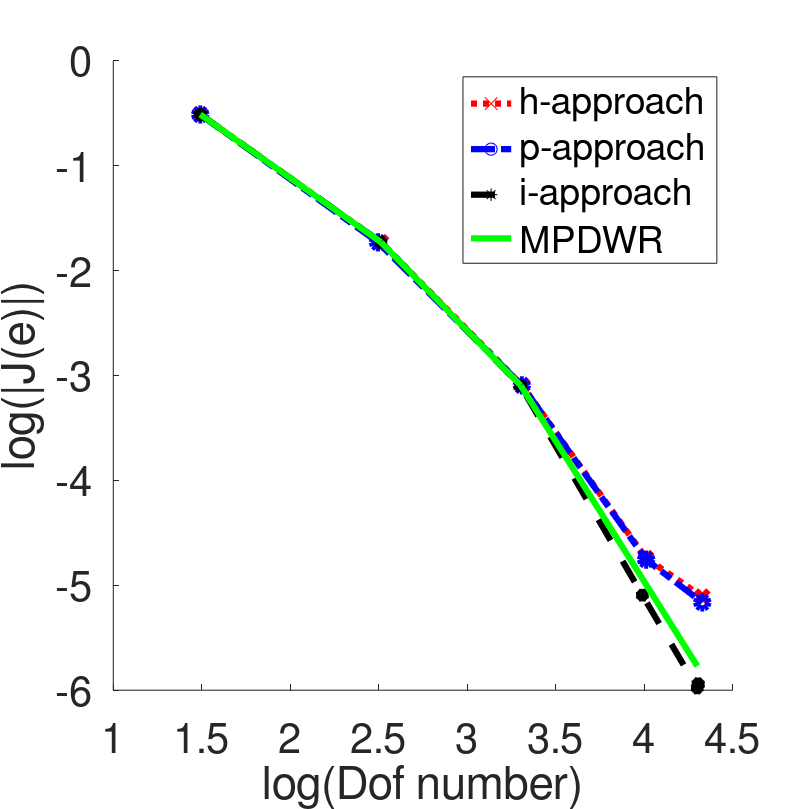}
\end{minipage}%
	\caption{Upper left: the adapted mesh based on MP-DWR. Lower left: the adapted mesh based on residual-based method. Middle: the relationship between error in target functionals and degrees of freedom. Right: the target error over the number of Dofs.}
\label{compare_result1}
\end{figure}

And the residual-based method with indicator $\eta_{K}^{resi}$ as mentioned above is implemented for comparison. The refined mesh and numerical results are shown in Figure \ref{compare_result1}.

From the left two pictures in Figure \ref{compare_result1}, the adapted mesh based on residual based error indicator gets denser around the singularity of function $u$ in the whole domain. By contrast, the adapted mesh based on the MP-DWR method not only gets denser around the singularity but also adds much more degrees of freedom around the interested point.

Moreover, the middle picture in Figure \ref{compare_result1} shows that the MP-DWR can obtain the numerical result reaching at the specified level with much more fewer degrees of freedom compared to the residual-based method. These phenomena confirm the goal-orientation property of the MP-DWR method, which demonstrates the effectiveness of the MP-DWR method on the Poisson problem with some particular target functionals somehow.

Furthermore, we also implement the $h$-approach, $p$-approach and $i$-approach for comparisons. The numerical results are shown in the right picture of Figure \ref{compare_result1}. It can be observed that the MP-DWR method can obtain a result of similar accuracy with a similar number of Dofs compared to the three classic approaches. Consequently, the efficiency of the MP-DWR method is comparable with the three approaches. However, the implementation of the MP-DWR method is faster, which makes the MP-DWR competitive. Such an advantage of the MP-DWR method will be shown in the following numerical example \ref{acceleration}.

\subsubsection{Coefficient elliptic equation}
Then we consider a coefficient elliptic equation, which is given by
 
\begin{equation}
	\left\{	
	\begin{array}{rl}
		&-\nabla\cdot (a\nabla u_3) + c u_3 = f,\text{ in } \Omega,\\
		&u_3\vert_{\partial\Omega} = 0,
	\end{array}
	\right.
\end{equation}
 
where
 
\begin{equation}
	a=
	\left\{ 
	\begin{array}{rl}
		&10,\ x_1^2-x_2^2\geq 0,\\ 
		& 1,\ x_1^2-x_2^2<0.
	\end{array}
	\right.
\end{equation}
 
And we set  $c = 1$ and $f = 1$. 
The output functional of this case is defined by
 
\begin{equation}
	J_p(u):=u(\boldsymbol{x}_e),\ \boldsymbol{x}_e=(0.5,0).
\end{equation}
 
\begin{figure}
	\begin{minipage}[c]{0.4\linewidth}
		\includegraphics[width=0.6\textwidth]{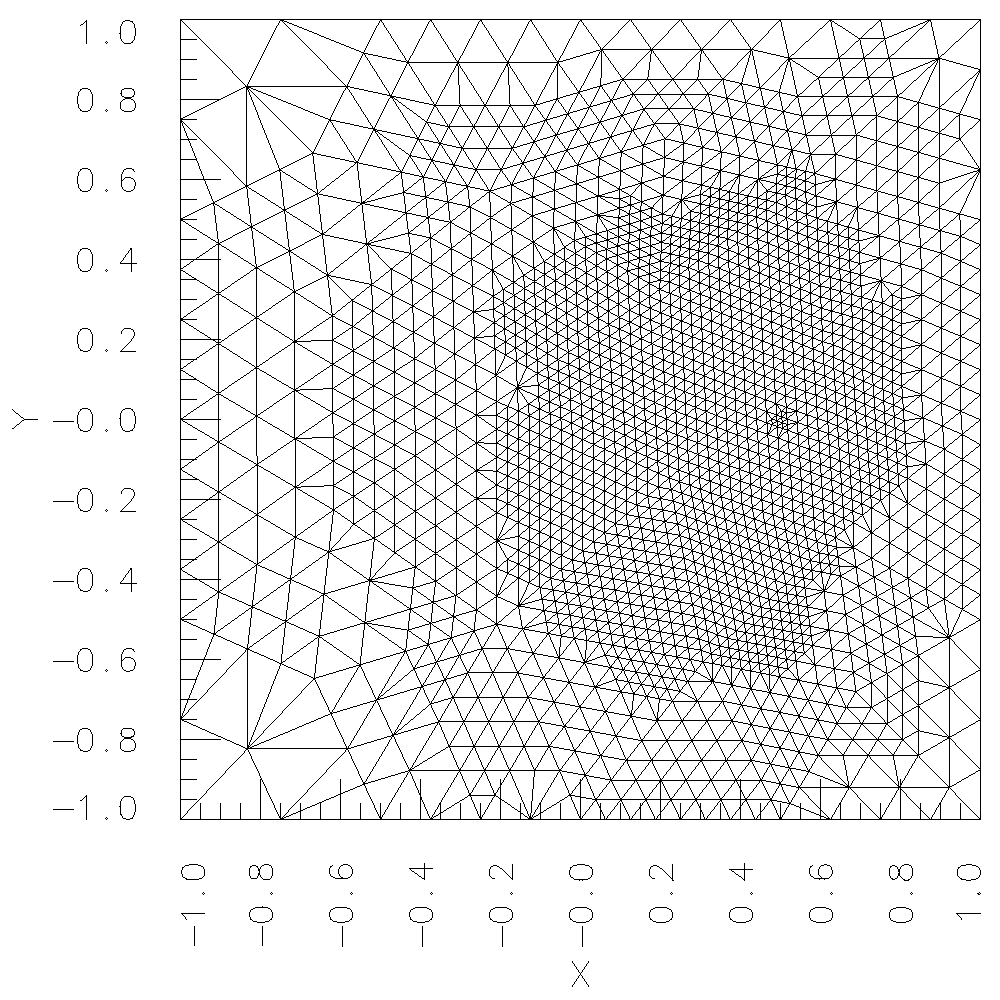}\\
		\includegraphics[width=0.6\textwidth]{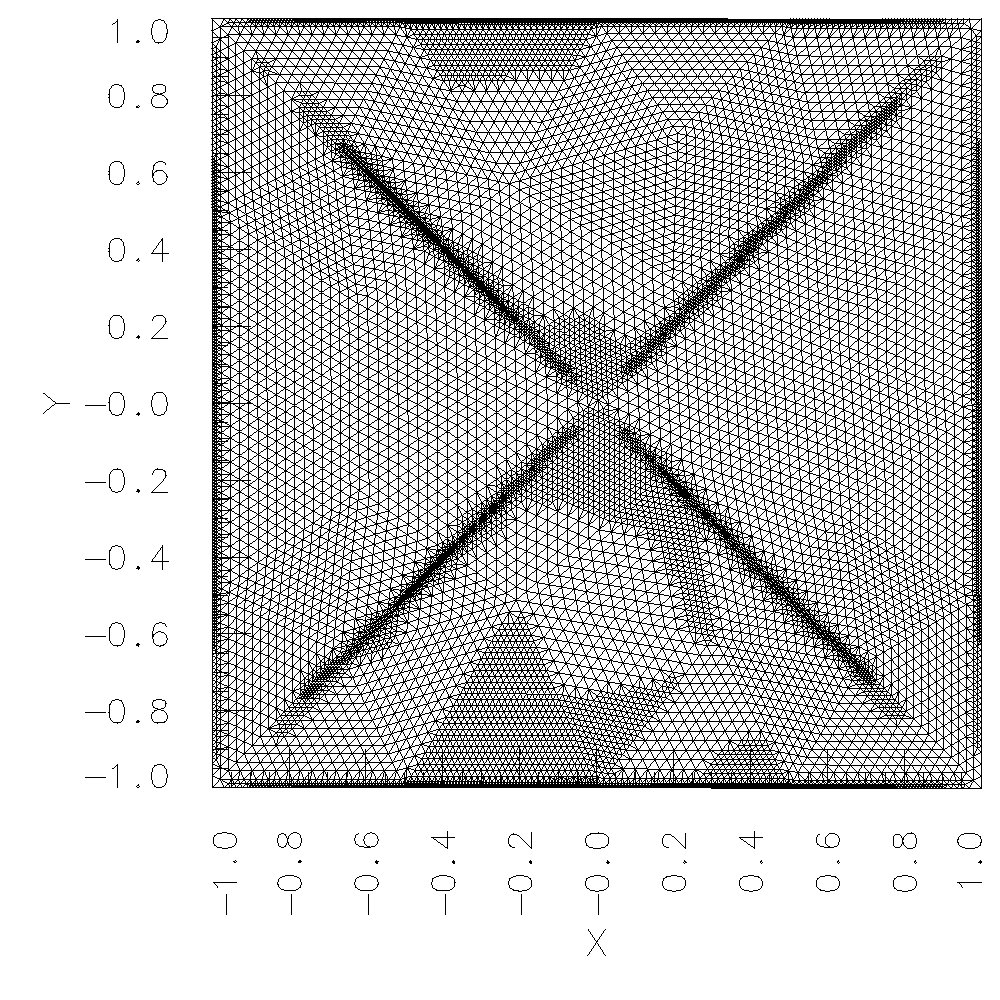}
	\end{minipage}
	\begin{minipage}[c]{0.6\linewidth}
		\includegraphics[width=0.7\textwidth]{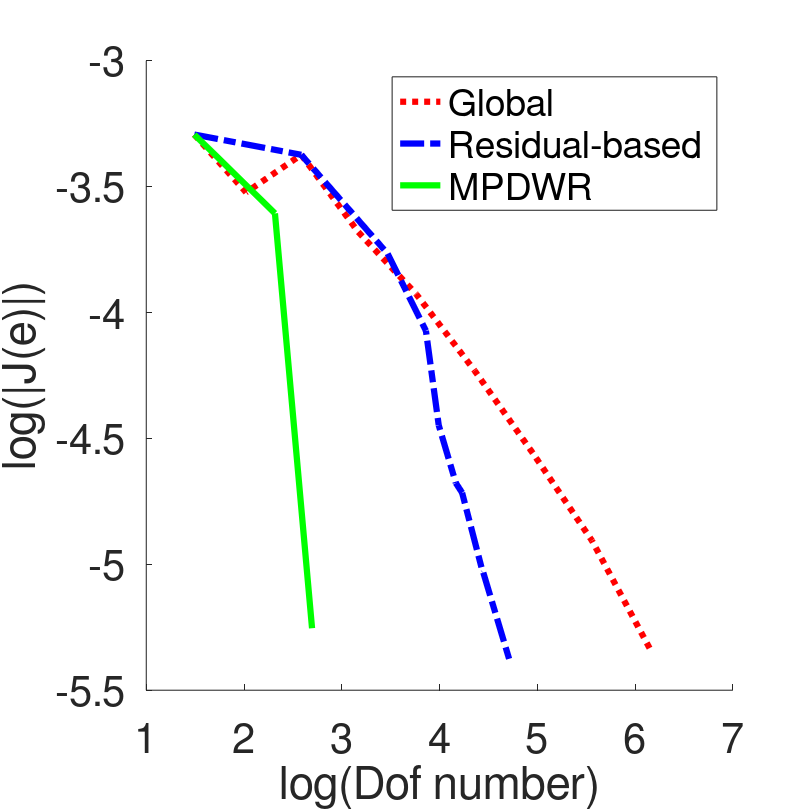}
	\end{minipage}%
	\caption{Upper left: the adapted mesh based on MP-DWR. Lower left: the adapted mesh based on residual-based method. Right: the log value of target functionals.}
	\label{compare_result_Coeff_Point}
\end{figure}

The bilinear form of this equation should also be changed. More details can be found in \cite{li2005multi}. Based on that, an error indicator similar to (\ref{dwrError}) can be obtained. We implement the MP-DWR method and residual-based method to refine the mesh grid, respectively. The right picture in Figure \ref{compare_result_Coeff_Point} shows that target error over the number of Dofs. Since the mesh is too coarse in the beginning, slight oscillation can be observed. However, with the refinement of the mesh grid, the convergence becomes smooth. Then it can be observed that the MP-DWR method can obtain the numerical result reaching the specified level with much fewer Dofs than the residual-based method. Besides that, from two left pictures in Figure \ref{compare_result_Coeff_Point}, we can observe that the refined mesh based on MP-DWR can capture the goal point better. While the residual-based method provides a mesh, on which the solution can be approximated well in the whole domain. As a result, it has to cost more expensive numerical effort than the MP-DWR method to get accurate quantities of interest. Based on these, the effectiveness of the MP-DWR method in the coefficient elliptic equation can be verified.
\subsubsection{Convection-diffusion equation}
To test the validity of MP-DWR method in other complex cases, we consider the convection-diffusion equation, which is given by
 
\begin{equation}
	\label{CovDiff}
	\left\{	
	\begin{array}{rl}
		-\nabla\cdot(\epsilon\nabla u) + \boldsymbol{\alpha}\nabla u = f,&\text{ in } \Omega,\\
		u = 0,&\text{ on } \partial \Omega,
	\end{array}
	\right.
\end{equation}
 
and the exact solution $u$ is given as $u_4(x_1,x_2)=sin(2\pi x_1)sin(\pi x_2)$.
In the numerical simulation, the $\epsilon$ is set as 0.01 and $\boldsymbol{\alpha} = [1,1]$. Following the DWR framework in the convection-diffusion equation, we can get the corresponding DWR error indicator for mesh adaptation. More details can be found in \cite{DWRBook,bruchhauser2017numerical}. For comparison, residual-based method with $\eta_K^{resi}$ is implemented in the simulation. 
In this case, the target functional is defined as the point-wise one:
 
\begin{equation}
	J_p(u):=u(\boldsymbol{x}_e),\ \boldsymbol{x}_e=(-\frac{1}{3},-\frac{1}{2}).
\end{equation}
 
\begin{figure}
	\begin{minipage}[c]{0.4\linewidth}
		\includegraphics[width=0.6\textwidth]{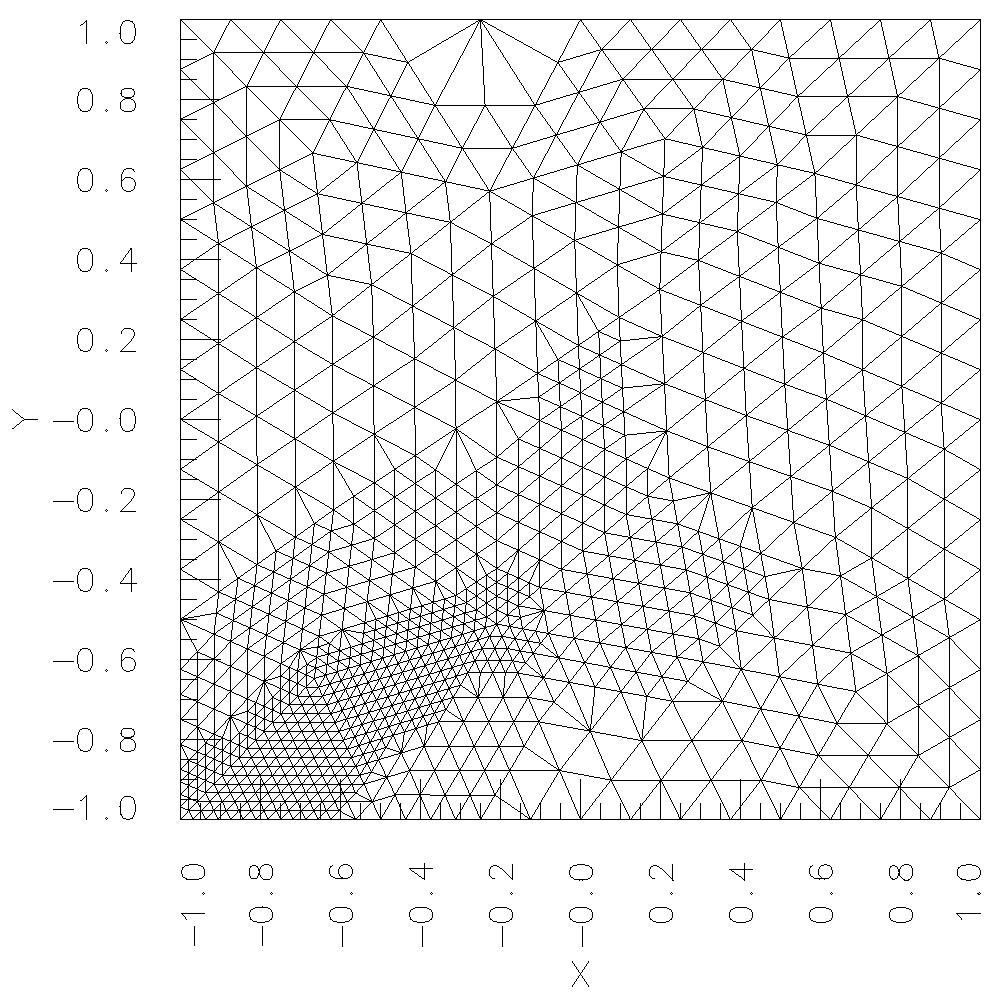}\\
		\includegraphics[width=0.6\textwidth]{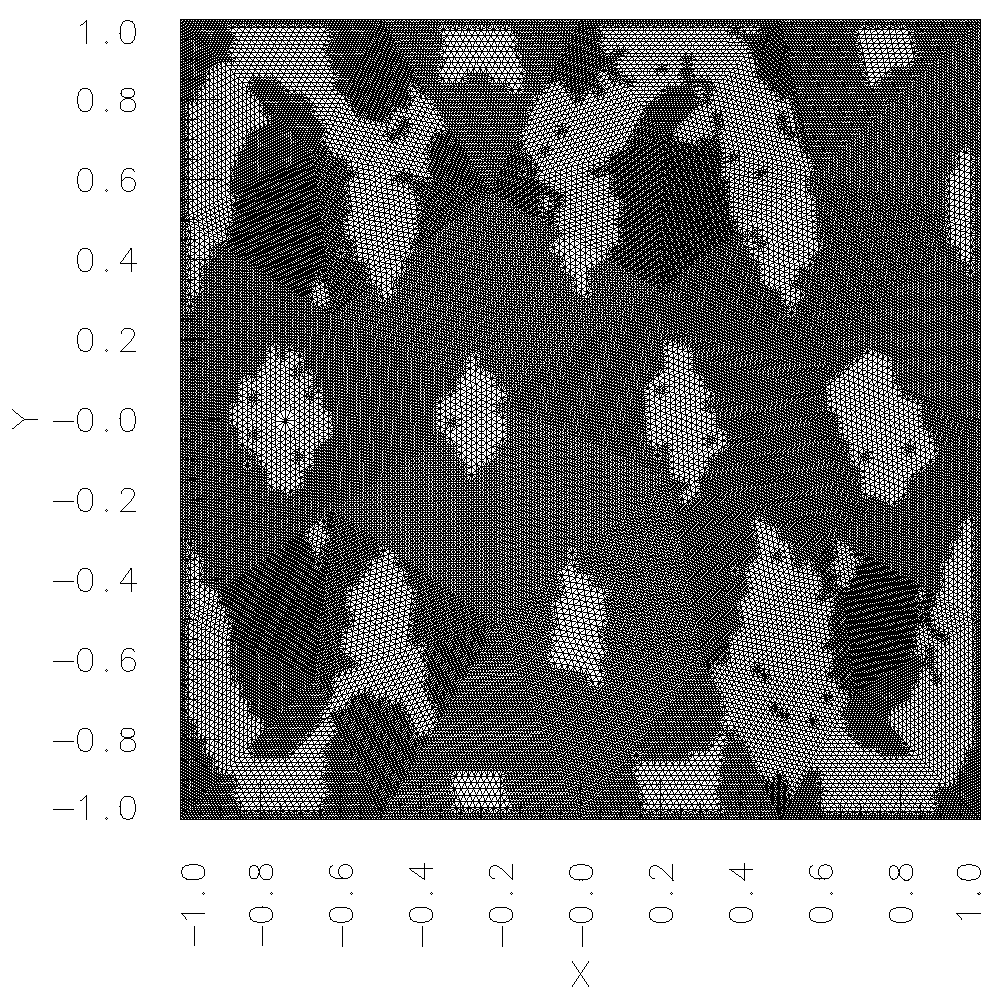}
	\end{minipage}
	\begin{minipage}[c]{0.6\linewidth}
		\includegraphics[width=0.7\textwidth]{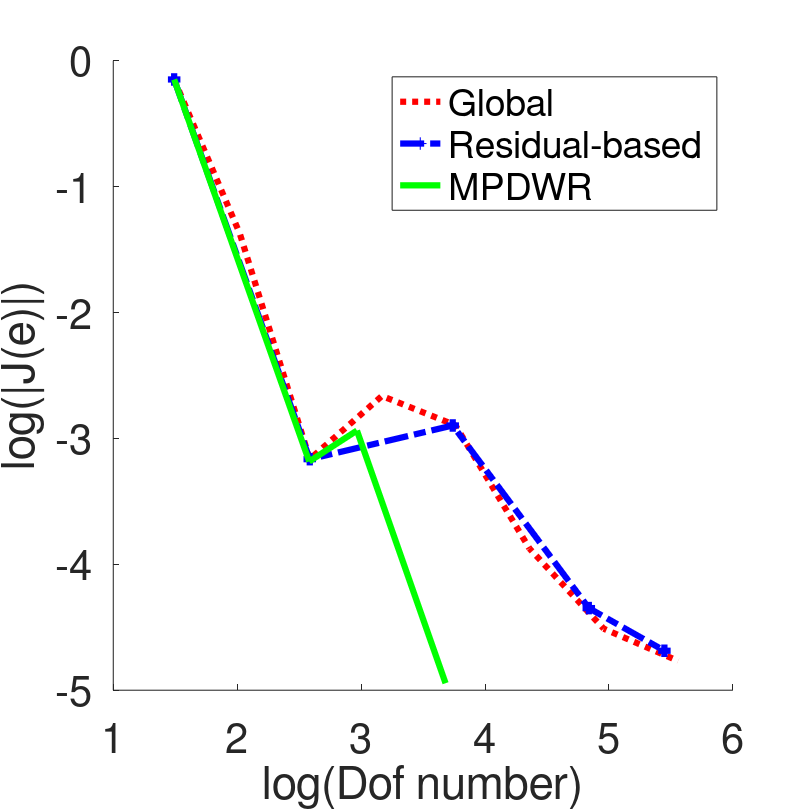}
	\end{minipage}%
	\caption{Upper left: the adapted mesh based on the MP-DWR method. Lower left: the adapted mesh based on the residual-based method. Right: the relationship between error in target functionals and degrees of freedom.}
\label{CovPoint_result}
\end{figure}
Numerical results are shown in the Figure \ref{CovPoint_result}. This case nicely shows the goal-orientation property of MP-DWR method. It can be seen that the MP-DWR method sets more degrees of freedom around the interested point and along those finite elements that influence the error of interested point in the direction of the flow field $\boldsymbol{\alpha}$. As a comparison, the residual-based method puts more degrees of freedom on the whole domain. Consequently, the adapted mesh grids would be quite dense as the lower left picture in Figure \ref{CovPoint_result}. And the right picture in Figure \ref{CovPoint_result} demonstrates that the residual-based method has to use many degrees of freedom to get an accurate numerical result. While the MP-DWR method can obtain the numerical result satisfying the same precision with just a few degrees of freedom. 
Therefore, we can confirm the effectiveness of MP-DWR method in this convection-diffusion case.

\subsection{Performance on acceleration}\label{acceleration}
\counterwithin{figure}{subsection}
\counterwithin{table}{subsection}
In this part, we will show one important advantage of the MP-DWR method in numerical simulations, i.e. the acceleration of simulations. Three classic DWR approaches, i.e. $h$-approach, $p$-approach, and $i$-approach, are implemented during the experiments as comparisons. By comparing the CPU time consumption of each approach, we can test the acceleration capability of the MP-DWR method.

\subsubsection{2D case}
\counterwithin{figure}{subsection}
\counterwithin{table}{subsection}
Above all, we consider the Poisson case as (\ref{Poisson}) with point-wise target functional $J_p$ in (\ref{Target4Poisson}). The exact solution is $u=5.0\times(1-x_1^2)(1-x_2^2)exp(1-x_2^{-4})$ and the interested point is $\boldsymbol{x}_e=(0.2,0.8)$. In this case, the DWR error indicators based on four different approaches (MP-DWR and three classic approaches) are generated from the same initial mesh with almost thirty thousands degrees of freedom. The process is repeated 100 times and the average CPU time of several operations in this process is recorded in Table \ref{CPUtime4Poisson}. 

It can be seen that all approaches cost similar CPU time to obtain the primal solution. However, in order to generate the error indicator (including solving the dual problem and constructing the error indicator), four approaches cost quite different CPU time.
Compared with $h$- and $p$-approach, the $i$-approach and MP-DWR method can save CPU time by avoiding several operations of building new space. Moreover, since no more Dofs or mesh grid points are introduced, they cost much less time to build quadrature information and solve the linear system than the former two classic approaches. In addition, the $i$-approach has to spend some time on interpolation while constructing the error indicator. Due to such a difference, the MP-DWR method is slightly faster than the $i$-approach.
\begin{table}[h]
	\centering
	\sffamily
	\begin{tabular}{|l|l|l|l|l|}
		\hline
		Operation                                    & $h$-approach           & $p$-approach          & $i$-approach & MP-DWR       \\ \hline
		Build primal FEM space                       & 8.915e0 & 8.274e0& 7.825e0 & 8.154e0 \\ \hline
		Build Quadrature information                 & 2.632e1 & 2.658e1 & 2.255e1& 2.196e1  \\ \hline
		Solve Linear System                           & 1.788e2 & 1.789e2 & 1.836e2& 1.786e2 \\ \hline
		Remeshing                                    & 1.041e1 &             & &            \\ \hline
		Build dual FEM space                         & 3.476e1 &1.092e1 &     &        \\ \hline
		Assemble Linear System                 & 8.231e1 & 4.912e1 & 1.837e1& 1.832e1  \\ \hline
		Solve Linear System                           & 7.039e0 & 2.669e1 & 1.238e0& 1.203e0 \\ \hline
		Construct indicator                          & 2.999e1 & 2.934e1 & 3.420e1& 2.815e1 \\ \hline
		\textbf{Generate indicator} &  \textbf{1.645e2}  & \textbf{1.161e2} & \textbf{5.381e1}& \textbf{4.767e1}  \\ \hline
		Total CPU time &3.785e2 &3.298e2&2.678e2&2.564e2 \\ \hline
	\end{tabular}
	\caption{\small CPU time(seconds) of operations in the computation of Poisson case. The values are the average of 100 experimental results. The process of generating the indicator includes solving the dual problem and constructing the error indicator. For the $i$-approach, the CPU time for interpolation is included in the CPU time for the `construct indicator' step.}
	\label{CPUtime4Poisson}
\end{table}

\begin{figure}[h]
	\centering
	\includegraphics[width=0.4\textwidth]{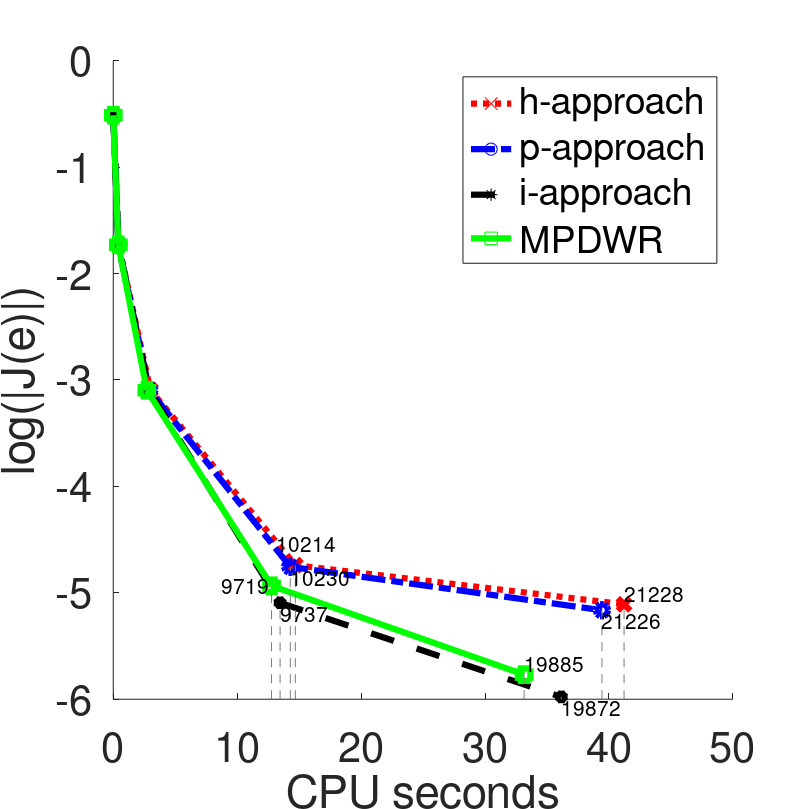}\includegraphics[width=0.4\textwidth]{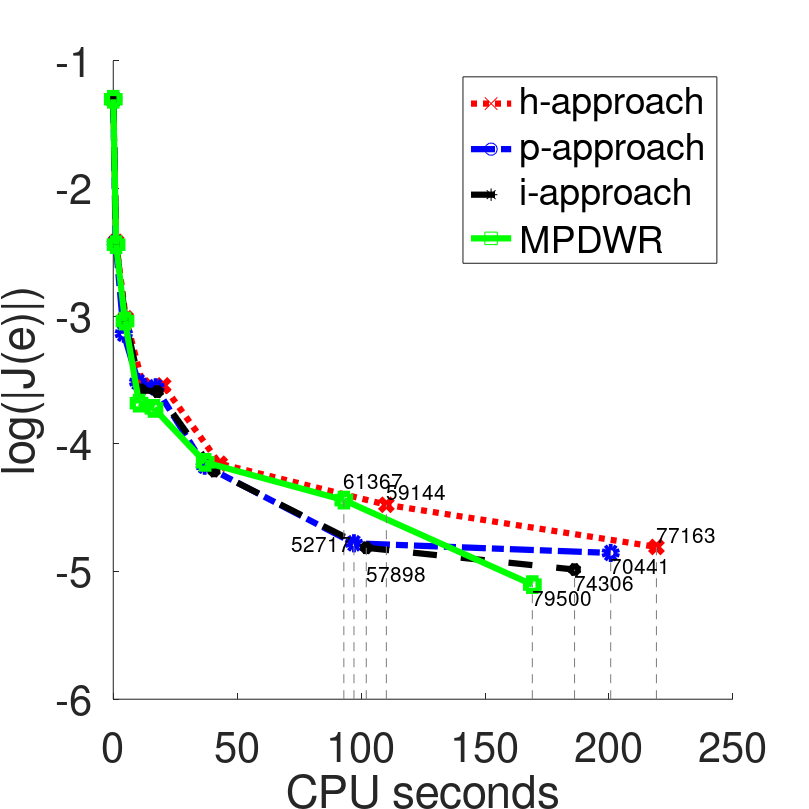}
	\caption{The error and CPU time of the whole experiments based on MP-DWR and three classic DWR approaches. Left: with point-wise target functional. Right: with area target functional. The number of Dofs is shown in the figure.}
	\label{FigCPU4Poisson}
\end{figure}

Then we implement the whole numerical simulation through each approach. The mesh adaptation is executed until the error tolerance is satisfied. The numerical results are shown in the left one of Figure \ref{FigCPU4Poisson}.
Furthermore, we change the exact solution of Poisson case as $u_1=(1-x_1^2)(1-x_2^2)\sin(4x_1)\sin(4x_2)$ with the area target functional $J_{\text{area}}(u)$ in (\ref{Target4Poisson}). Similar simulations are implemented. 
The relationship between the error and whole CPU time for this case is shown in the right one of Figure \ref{FigCPU4Poisson}.
In the left picture in Figure \ref{FigCPU4Poisson}, it can be seen that the MP-DWR method and $i$-approach are much faster than $h$- and $p$-approach. The acceleration gets more obvious as the times of mesh adaptation increase. Then we focus on the comparison of MP-DWR and $i$-approach: the $i$-approach can control the error of the point-wise target functional better than MP-DWR method with a similar number of Dofs. But the CPU time of the MP-DWR method is always less than the $i$-approach. This suggests that the MP-DWR method may be less accurate but computationally fast compared to the $i$-method in this case. 

As for the right picture, the corresponding points of the MP-DWR method always stay at the most left position in each iteration step, which means that the MP-DWR method generally cost the least CPU time in simulations. When the iterations stop, although the MP-DWR method uses the most Dofs, the CPU time of the MP-DWR method is significantly less than the other three methods. This demonstrates that the MP-DWR method can be calculated faster than the three classical methods in this case.

\subsubsection{3D case}
\counterwithin{figure}{subsection}
\counterwithin{table}{subsection}
Then we consider the Poisson equation in 3-dimensional with a target functional, which is given by:
\begin{equation}
	J_{G}(u)=\int_{\Omega_{3D}} u d\boldsymbol{x},\ \boldsymbol{x}=(x_1,x_2,x_3),
\end{equation}
where $\Omega_{3D}=(-1,1)^{3}$.
The exact solution is given as the example in \cite{DWRBook}:
\begin{equation}
u=(1-x_1^2)^2(1-x_2^2)^2(1-x_3^2)^2 (64x_1^2+0.1)^{-1} (64x_3^2+0.1)^{-1}.
\end{equation}
The comparison of CPU time for the whole simulation based on the MP-DWR method and three classic approaches is similar to the 2D case in Figure \ref{FigCPU4Poisson}. The MP-DWR method and $i$-approach demonstrate similar efficiency, and both of them are computed faster than the $h$- and $p$-approach. 

Then we mainly focus on the comparison of CPU time for generating the error indicator. Similar to the 2D case, we implement the four approaches to generate the error indicators from the same initial mesh (with almost 18000 Dofs), respectively. The corresponding CPU time is recorded in the Table \ref{3DGenerate}. The values are the average of 100 experimental results.

Compared to the $i$-approach, the MP-DWR method can save some CPU time during solving the dual problem since the involvement of single-precision computations. Besides that, due to the increased number of elements and corresponding neighbors in the 3D case, higher-order interpolation takes more CPU time than in the 2D case. Consequently, the MP-DWR method can save much CPU time during constructing the indicator. More importantly, compared to the $h$- and $p$-approach, the improvement in CPU time brought by the MP-DWR method is quite significant. The MP-DWR method can use the equivalent of about 15\% of the time for the $h$-approach and about 30\% of the time for the $p$-approach to generate the indicator. Based on that, we can demonstrate that the performance of the MP-DWR method in acceleration is quite considerable in the 3D case compared to the $h$- and $p$-approach. Therefore, the MP-DWR method can be more competitive in 3D numerical simulations.

\begin{table}[h]
	\centering
	\sffamily
	\begin{tabular}{|l|l|l|l|l|}
		\hline
		Operation                                    & $h$-approach           & $p$-approach           & $i$-approach& MP-DWR       \\ \hline
		Build primal FEM space                       & 3.086e0 & 3.048e0 & 2.825e0& 3.196e0 \\ \hline
		Build Quadrature information                 & 1.135e1 & 1.137e1 & 1.129e1& 9.710e0 \\ \hline
		Solve Linear System                           & 2.137e1 & 2.137e1 & 2.139e1& 2.148e1 \\ \hline
		Remeshing                                    & 6.423e0 &             &          &   \\ \hline
		Build dual FEM space                         & 2.446e1 & 3.886e0 &   &          \\ \hline
		Assemble Linear System                 & 7.833e1 & 4.096e1  & 9.853e0& 8.272e0\\ \hline
		Solve Linear System                           & 1.062e0 & 6.402e0 & 1.449e-1& 1.447e-1 \\ \hline
		Construct indicator                          & 9.956e0 & 1.033e1 & 2.157e1 & 9.691e0  \\ \hline
		\textbf{Generate indicator}           & \textbf{1.202e2 }& \textbf{6.162e1} &  \textbf{3.156e1}& \textbf{1.811e1}  \\ \hline
		Total time							         & 1.560e2 & 9.741e1 & 6.707e1& 5.249e1  \\ \hline
	\end{tabular}
	\caption{\small CPU time(seconds) of operations in 3-dimensional Poisson case. The values are the average of 100 experimental results. The process of generating the indicator includes solving the dual problem and constructing the error indicator. For the $i$-approach, the CPU time for interpolation is included in the CPU time for the `construct indicator' step.}
	\label{3DGenerate}
\end{table}

\subsection{Performance in eigenvalue problems}
\begin{figure}
	\begin{minipage}[c]{0.45\linewidth}
		\includegraphics[width=0.5\textwidth]{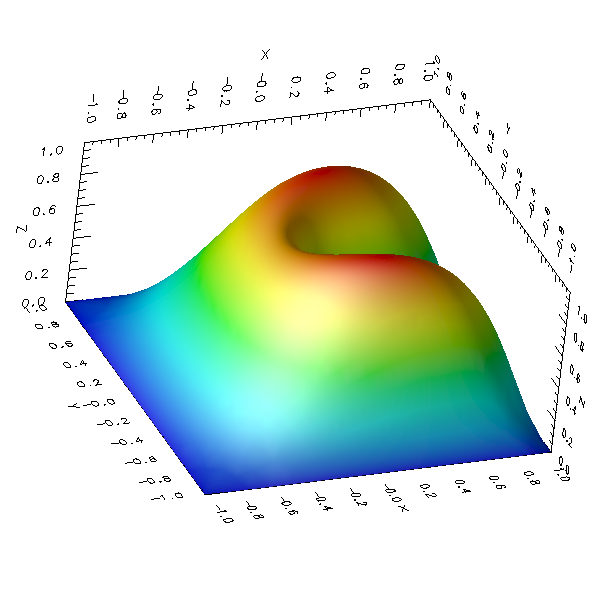}\\
		\includegraphics[width=0.5\textwidth]{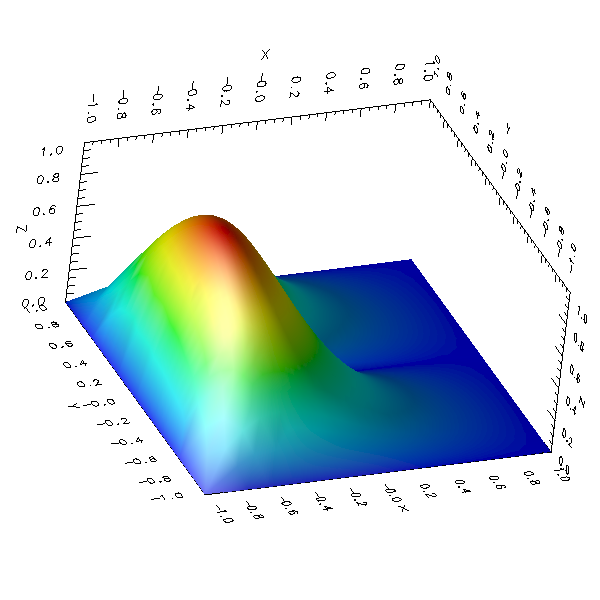}
	\end{minipage}
	\begin{minipage}[c]{0.6\linewidth}
		\includegraphics[width=0.6\textwidth]{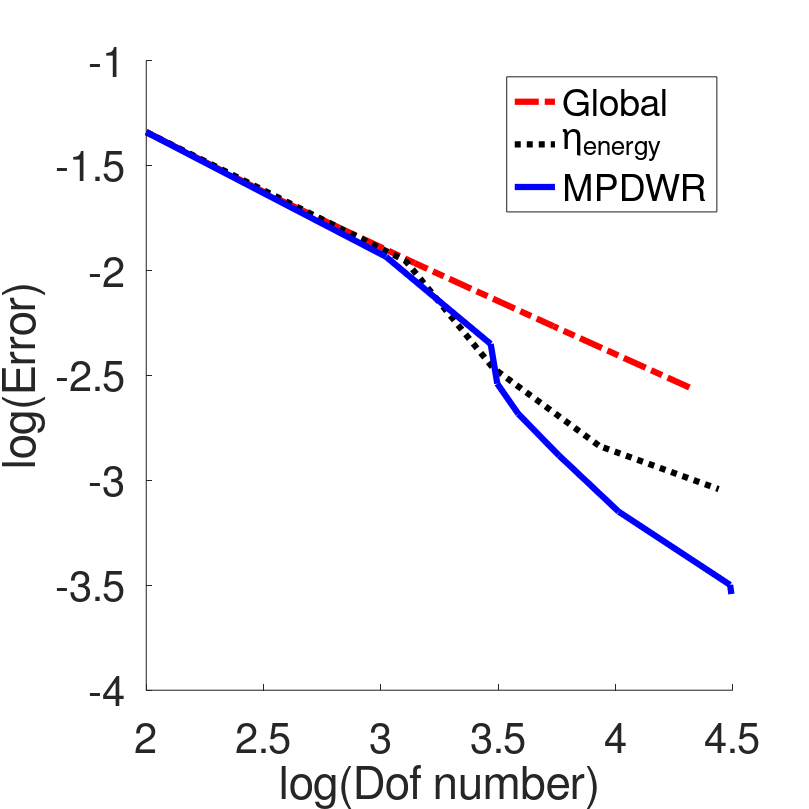}
	\end{minipage}%
	\caption{Upper left: the plot figure of the primal solution. Lower left: dual solution.. Right: the smallest eigenvalue errors over the number of Dofs.}
	\label{EigEffectivity}
\end{figure}

\counterwithin{figure}{subsection}
\counterwithin{table}{subsection}
To extend the MP-DWR method into more applications, we spend some time on verifying its effectiveness in the eigenvalue problem. Consider the eigenvalue problem as follows. Seeking an eigenpair
$(\lambda,u)\in\mathcal{C}\times H_{0}^{1}(\Omega)$ with
\begin{align}
	\label{eig_problem}
	 -\Delta u + \boldsymbol{\beta} \cdot \nabla u= \lambda u\ in\ \Omega.
\end{align} 
In experiments, $\boldsymbol{\beta}$ is set as [3,0]. The computational domain is  set as the slit area $\Omega_3=((-1,1)\times(-1,1))\verb|\|([0,1]\times\{0\})$. Following the DWR framework in the eigenvalue problem as in \cite{DWRBook,heuveline2001posteriori,Gedicke2013On}, we can get the related error indicator for mesh adaptation. In the simulation of the MP-DWR method, we follow the framework of Algorithm \ref{alg:Framwork}, i.e. the primal solution and corresponding residual $\rho_{K}$ are computed in double-precision. While the dual solution and corresponding residual $\rho_{K}^{*}$ are computed in single-precision. Based on these, the error indicator $\eta_{\lambda}^{\omega}=\sum_{K\in\mathcal{T}_h}\{\rho_{K}\omega_{K}^{*}+\rho^{*}_K\omega_{K}\}$ can be computed. For comparison, the residual-based error indicator $\eta_{\lambda}^{\text{energy}}=c_{\lambda}\sum_{K\in\mathcal{T}_h}h^{2}_{K}\{\rho_{K}^2+\rho^{*2}_K\}$ is also implemented, where $c_\lambda$ is a constant growing linearly with $\vert\lambda\vert$. The left pictures in Figure \ref{EigEffectivity} shows the numerical primal and dual solution. The residual of the smallest eigenvalue $\lambda_1$ over Dofs is shown in the right one of Figure \ref{EigEffectivity}. It can be observed that, with the specified number of Dofs, the MP-DWR method can get a more accurate result than the $\eta^{\text{energy}}$ based method. This demonstrates that the MP-DWR method can deal with the interested eigenpair more efficiently. Consequently, it can be demonstrated that the MP-DWR method can still work in the eigenvalue problem. Based on these, it can be expected that the MP-DWR method can be used in many applications, such as density functional theory.

\subsection{The limit of MP-DWR method}\label{limitofFD}
\begin{figure}[!htbp]
	\centering
	\includegraphics[width=0.6\textwidth]{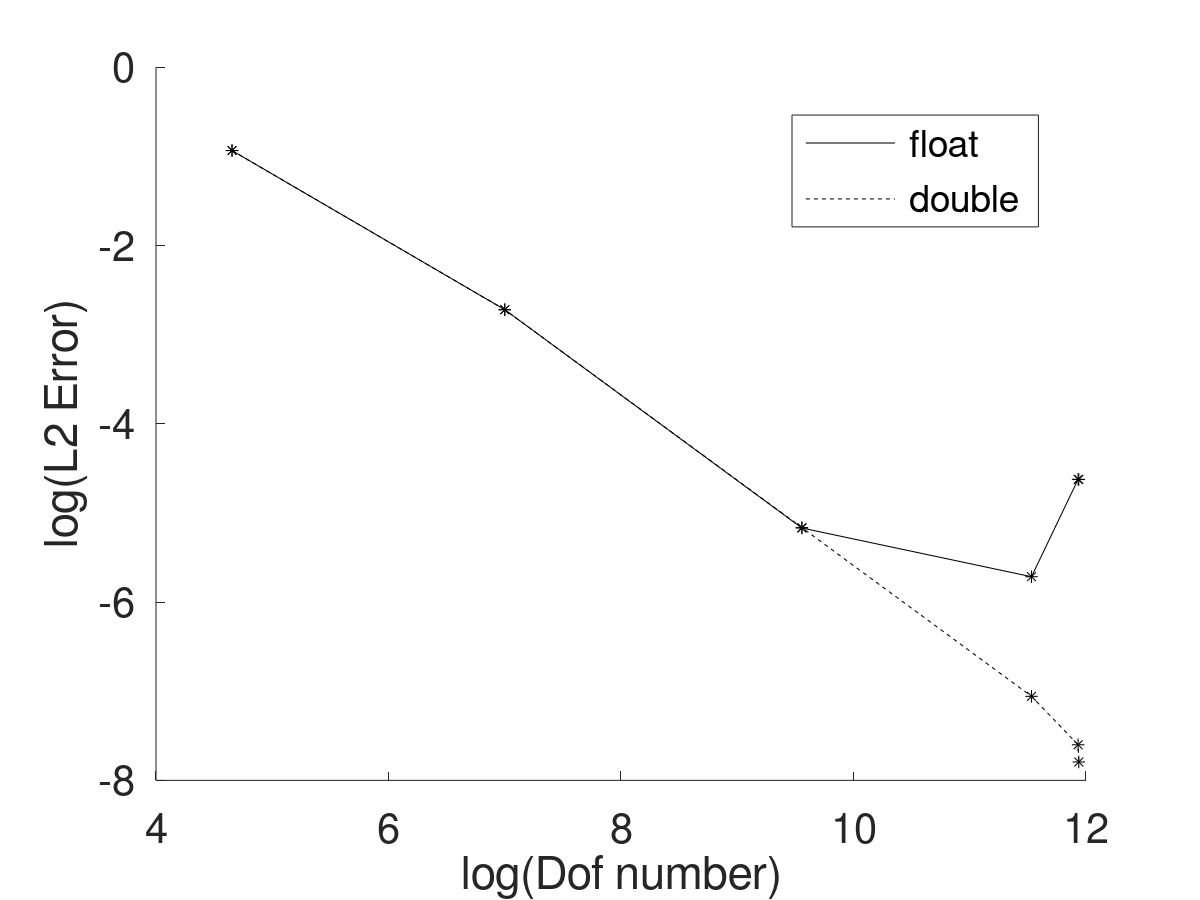}
	\caption{Primal solutions computed in float and double precision.}
	\label{limit}
\end{figure}

In this part, we discuss a limit of the MP-DWR(FD) framework. In some numerical experiments, it is strange that MP-DWR(FD) becomes less effective
with refining the mesh too many times. Consider the Poisson case with point-wise target functional in case \ref{PoissonCase}. We solve the primal
problem with single-precision and double-precision, respectively.  The numerical results
of the primal solutions are shown in
Figure \ref{limit}. It can be seen that when the mesh grids are dense enough,
the residual in single-precision gets larger as the number of degrees of freedom increases, which is different from the residual in double-precision. With this unexpected phenomenon, it can not be expected that such an inaccurate primal solution in single-precision can construct a reliable DWR error indicator.
One possible reason is that while we adapt the mesh several times, some area will gets quite dense and the minimal volume of elements will get
very tiny. For reference, the minimal volumes of elements on the mesh grids in
each iteration are shown in Table \ref{volume table}. When the minimal volume
gets the level of 1.0e-6, the data in such an element can be close to the limit of single-precision so that the computations may demonstrate some unexpected phenomena causing the strange curve in Figure \ref{limit}.
Due to such a limit, we have to be careful about the data precision in numerical simulations while using the MP-DWR method. 
\begin{table}[h]
	\begin{center}
		\label{volume table}
		\begin{tabular}{@{}llllll@{}}
			\toprule
			Adapt times	& 1 & 2 & 3 & 4 & 5\\ 
			Minimal element volume  & 1.617e-2 & 1.010e-3 & 6.316e-05 & 4.34e-06 & 1.247e-06\\ 
			\bottomrule
		\end{tabular}
			\caption{The minimal volume of elements in adaptation}
	\end{center}
\end{table}

\begin{remark}
  Although there is a limit to the MP-DWR(FD) method, we propose a potential
  algorithm to make full use of the advantages of the MP-DWR method: First of all,
  we can compute the primal solution with half-precision and dual solution with
  single-precision, respectively. When the primal solution gets the limit of
  half-precision, we change the precision parameters in our programming to
  compute the primal solution with single-precision and dual solution with
  double-precision, respectively. In this way, we can make full use of the
  acceleration of the MP-DWR method to get a reliable numerical result. And such an algorithm can also be an efficient preconditioned approach.
\end{remark}

\section{Conclusion}
To accelerate the process of obtaining dual solution in the DWR methods, we
utilize the properties of precision formats and propose a novel DWR
implementation based on multiple-precision named MP-DWR. Numerical experiments
confirm the effectiveness of the novel approach. Memory storage savings can be expected due to the low-precision calculations involved. Furthermore, considerable reduction in CPU time compared with classic $h$- and $p$-approach is shown in simulations, which makes the MP-DWR method quite competitive. Compared with $i$-approach, the MP-DWR may lose some accuracy but it gains some acceleration in simulations, especially in the 3D case. By combining the MP-DWR method and $i$-approach, both of the accuracy and efficiency might be improved in simulations, which is our ongoing work. Although there is a limit to the MP-DWR method due to the property of precision formats in some situations, the substantial improvements in data storage and CPU time suggest that it is still an efficient method. 

In addition, some more precise floating-point data type than double, such as \_\_float128(long double), may be applied in the MP-DWR framework to overcome the limit. Such an extension of the MP-DWR method will be implemented in our future study. In our previous work, applications of the $h$-adaptive method to computational fluid dynamics, micromagnetics, and density functional theory have been studied \cite{BAO20124967,bao2013numerical,hu2016adjoint,sci2019AnAF,Meng2022ANF}. The introduction of the MP-DWR method into these applications to further improve the efficiency is also our future work. Besides that, the dual weighted residual method has been applied to many time-dependent problems \cite{kocher2019efficient,bruchhauser2020dual,SLEEMAN2022114206}, for which the MP-DWR is expected to be used to improve the computational efficiency.

\section*{Acknowledgement}
The research of G. Hu was partially supported by National Natural Science
Foundation of China (Grant Nos. 11922120 and 11871489), FDCT of Macao SAR
(0082/2020/A2), MYRG of University of Macau (MYRG2020-00265-FST) and
Guangdong-Hong Kong-Macao Joint Laboratory for Data-Driven Fluid Mechanics and
Engineering Applications (2020B1212030001).

\section*{Statements and Declarations}
\subsection*{Data availability statements}
Data sharing not applicable to this article as no datasets were generated or 
analysed during the current study.

\subsection*{Competing interests}
The authors declare that they have no conflict of interest.

\bibliographystyle{plain}
\bibliography{bibliography}


\end{document}